\documentclass[a4paper]{article}
\usepackage[latin1]{inputenc}
\usepackage[english,french]{babel}
\usepackage[latin1]{inputenc}
\usepackage[T1]{fontenc}
\usepackage{graphicx}
\usepackage{amssymb}
\usepackage{amsthm}
\usepackage{amsmath,amssymb}
\usepackage[all]{xy}
\usepackage{enumitem}

\title{ D Théorème  D d'ouvert de bifurcations pour des automorphismes de Hénon 
 eT intersections d'ensembles de Cantor}
\author{Sébastien Biebler}
\date{Decembre 2013}
\AtEndDocument{    \footnotesize \textsc{LAMA, UMR8050, UNIVERSITE PARIS-EST MARNE LA VALLEE, 5 BOULEVARD DESCARTES, 77454 CHAMPS SUR MARNE, FRANCE} \par \textit{E-mail address} : \texttt{sebastien.biebler@u-pem.fr}}
\usepackage[nottoc, notlof, notlot]{tocbibind}
\begin{document}
 \newtheorem* {myTheo2} {Theorem}  \newtheorem* {myTheo33} {Theorem A}
\newtheorem* {myTheo35} {Theorem B}
\newtheorem* {myTheo37} {Corollary C}

 \newtheorem {df5} [subsubsection]{Definition}  \newtheorem {nt5} [subsubsection] {Notation} \newtheorem {co5}[subsubsection]{Corollary} \newtheorem* {rk}{Remark}  \newtheorem* {coo}{Corollary} 
\newtheorem{df57}[subsection]{Definition} \newtheorem*{df557}{Definition}
\newtheorem{pr57}[subsection]{Proposition}
\newtheorem{lem3}[subsection]{Lemma}
\newtheorem{nt57}[subsection]{Notation}
\newtheorem{lem57}[subsection]{Lemma}
\newtheorem{rem57}[subsection]{Remark}

\begin{center} \Huge{A Complex Gap Lemma}
\end{center} \begin{center}  \huge{S\'ebastien Biebler } \end{center}
\selectlanguage{english}
\begin{abstract}   
Inspired by the work of Newhouse in one real variable, we introduce a relevant notion of thickness for dynamical Cantor sets of the plane associated to a holomorphic IFS. Our main result is a complex version of Newhouse's Gap Lemma : we show that under some assumptions, if the product $\mathrm{t}(K)\mathrm{t}(L)$ of the thicknesses of two Cantor sets $K$ and $L$ is larger than 1, then $K$ and $L$ have non empty intersection. Since in addition this thickness varies continuously, this gives a criterion to get a robust intersection between two Cantor sets in the plane.

 \end{abstract}
\tableofcontents 
\section{Introduction and main results} \subsection{Historical context} We say that a smooth dynamical system $f : M \rightarrow M$ defined on a manifold $M$ is uniformly hyperbolic if periodic points of $f$ are dense in the non-wandering set $\Omega(f)$ and if the tangent space along $\Omega(f)$ admits a continuous invariant splitting $E_{s} \bigoplus E_{u}$ where vectors in $E_{s}$ (resp. $E_{u}$) are uniformly contracted under $f$ (resp. under $f^{-1})$. The dynamics of such systems is particularly well understood (see \cite{caga}). \medskip

It was initially believed in the 1960's that uniform hyperbolicity was a dense property in the space of diffeomorphisms of a manifold. The discovery in the 1970's (see \cite{n}) of the so-called Newhouse phenomenon, i.e. the existence of residual sets of $C^{2}$-diffeomorphisms of compact surfaces with infinitely many sinks (periodic attractors) showed that this expectation was false. \medskip

 The main point in Newhouse's proof is the creation of robust homoclinic tangencies. To produce them, Newhouse takes a diffeomorphism which has a horseshoe. Newhouse isolates a certain curve called line of tangency with the property that the stable and unstable laminations of the horseshoe intersect this curve along two respective Cantor sets and any intersection between these Cantor sets corresponds to a tangency between the laminations. In particular a robust intersection gives rise to a robust tangency. To get a robust intersection of Cantor sets, Newhouse uses a technical result called the Gap lemma. As a first step, Newhouse assigns to any one-dimensional Cantor set $K$ a number $\tau(K)$ called thickness which measures to what extent the Cantor set $K$ is fat :  

 \begin{df557}
    For every Cantor set $K \subset \mathbb{R}$, a gap of $K$ is a connected component $I$ of $\mathbb{R} \backslash K$. We denote $\tau_{I} = \frac {\ell(U)} {\ell(I)}$, where $U$ is the smallest interval between $I$ and a gap of length equal or larger than the length $\ell(I)$ of $I$. Then, we denote by $\tau(K)$ the thickness of $K$ : $$\tau(K) = \mathrm{inf}_{I}  \tau_{I} $$ 
    \end {df557} The thickness of a Cantor set is related to its Hausdorff dimension : $$\dim_{H}(K) \ge \frac{ \log 2}{\log(2+\tau(K)^{-1})}$$ In particular, Cantor sets with high thickness have a Hausdorff dimension close to 1. Then, a geometric proof (see \cite{bb5}) gives :

\begin {myTheo2}[Newhouse's Gap lemma] Let $K \subset \mathbb{R}$ and $L \subset \mathbb{R}$ be two one-dimensional Cantor sets such that : $\tau(K)\tau(L) \ge 1$. Then, one of the following is true : \begin{enumerate} \item or  $K$ is included in a gap of $L$ \item or $L$ is included in a gap of $K$ \item or  $K \cap L \neq \emptyset$ \end{enumerate}

                                                                              \end {myTheo2}

Since the thickness varies continuously and it is easy to get rid of the two first options, the previous result gives intersections which are in fact robust. Let us remark that other important  studies about intersections of Cantor sets in one dimension include \cite{31615} and \cite{316}. \medskip 

  In the complex setting, an extension of the Newhouse phenomenon was studied in the 1990's by Buzzard. Let us denote by $\text{Aut}_{d}(\mathbb{C}^{k})$ the space of polynomial automorphisms of $\mathbb{C}^{k}$ of degree $d$ for $d,k \ge 2$. Buzzard proved in \cite{bb1} that there exists $d > 0$ such that there exists an open set $N \subset \text{Aut}_{d}(\mathbb{C}^{2})$ such that automorphisms in $N$ have persistent homoclinic tangencies. The method is similar to the one-dimensional case, but this time one has to intersect two Cantor sets in the plane and not in the line. Buzzard gives an elegant criterion (see \cite{bb3}) which generates an intersection for two planar Cantor sets with very special geometric properties. Nevertheless, he does not give any general Gap lemma for planar Cantor sets, his construction is very specific and can't be applied for most Cantor sets. Let us also point out that Buzzard's work was extended  to higher complex dimensions in \cite{bieblerberthollove}, which relies on a different mechanism.  \medskip

\subsection{Main results}

In this article, we prove a complex Gap lemma for dynamical Cantor sets in $\mathbb{C}$ associated to Iterated Function Systems (IFS's) $\{f_{1},...,f_{p}\}$, where each $f_{i} : S \rightarrow f_{i}(S) \Subset S$ is defined on a square $S$ and univalent. To the dynamical Cantor set $K$ associated to such an IFS, we assign a thickness $\mathrm{t}(K)$ defined in Definition 2.14. We will work with dynamical Cantor sets satisfiying a condition called "well balanced" (see Definition 2.16). This condition is robust. It means that : \begin{enumerate} \item  the subsets $f_{i}(S)$ occupy a substantial part of the available space in the initial square $S$ \item all the subsets $f_{i}(S)$ are small compared to the initial square \end{enumerate}

The following Theorem is the main result in this paper and is a generalization of Newhouse's Gap lemma to the holomorphic context. It might be used to investigate the existence of the Newhouse phenomenon in spaces of polynomial automorphisms of given degree.

   \begin {myTheo33}
   Let $K$ and $L$ be two dynamical Cantor sets which are limit sets of respective IFS's defined on a same square. Suppose that $K$ and $L$ are well balanced and satisfy : 
$$\mathrm{t}(K)\mathrm{t}(L) \ge 1$$
Then $K \cap L$ is non empty.

   \end{myTheo33}
We will also show that the thickness is a continuous function of the maps defining the dynamical Cantor sets $K$ and $L$ :

 \begin {myTheo35}
 The thickness $\mathrm{t}(K)$ is a continuous function of the maps $\{f_{1}, \cdots, f_{p} \}$ defining $K$.

   \end{myTheo35}

These two results can produce a robust intersection of Cantor sets. Indeed, as an immediate corollary of Theorems A and B, we have : 

 \begin {myTheo37}
   Let $K$ and $L$ be two dynamical Cantor sets which are limit sets of respective IFS's defined on a same square. Suppose that $K$ and $L$ are well balanced and satisfy :
$$\mathrm{t}(K)\mathrm{t}(L) >1$$
Then $K$ and $L$ intersect in a robust way, this is if $K'$ and $L'$ are dynamical Cantor sets obtained after perturbation of the IFS's defining $K$ and $L$, then $K' \cap L' \neq \emptyset$.

   \end{myTheo37}

\subsection{Organization of the paper}

In Section 2, we define our notion of thickness for dynamical Cantor sets in $\mathbb{C}$ and we also define the "well balanced" condition. In Section 3, we give an example of two dynamical Cantor sets having a robust intersection. In Section 4, we prove Theorem A. Finally, we prove Theorem B in section 5. \newline \newline \textbf{Acknowledgments :}
The author would like to thank his PhD advisor, Romain Dujardin as well as Pierre Berger for useful comments. This research was partially supported by the ANR project LAMBDA, ANR-13-BS01-0002.

\section{Definition of the thickness}

In this article, we will work with dynamical Cantor sets in $\mathbb{C}$ associated to IFS's $\{f_{1},...,f_{p}\}$. In particular, we first need to define a norm on $\mathbb{C}$ and a suitable topology for the maps $f_{i}$. \medskip

\textbf{Choice of a norm on $\mathbb{C}$ : }
We use the usual euclidean norm $ | \cdot | = | \cdot |_{2}$ on $\mathbb{C}$ in all this article. For every $z = x+iy \in \mathbb{C}$, $|z| = |z|_{2} = \sqrt{x^{2}+y^{2}}$. All the distances will be measured relatively to this norm. \medskip 

 We will often use squares in the following. All these squares will be open. In particular, the diameter of a square $S$ is measured relatively to $|\cdot|$ and is equal to the length of its diagonal. A square $S$ of diameter 1 has its four sides of length $\frac{1}{\sqrt{2}}$ and contains an inscribed disk of diameter $\frac{1}{\sqrt{2}}$.  \medskip

We will also often use disks in the following. All these disks will be open. In particular, the diameter $\delta(\Gamma)$ of a disk $\Gamma$ is measured relatively to $|\cdot|$ and is equal to twice its radius.  \medskip

\textbf{Choice of a topology on the space of IFS's } : We will use IFS's $\{f_{1},...,f_{p}\}$ defined on a square $S \subset \mathbb{C}$, where every $f_{i}$ is a holomorphic contraction defined on $S$ such that $f_{i}(S) \Subset S$. We take the product topology of the compact open topology on the set of holomorphic maps from $S$ to $S$ to get a topology on the set of IFS's $\{f_{1},...,f_{p}\}$

 \begin{df57}
Let $\{f_{1},...,f_{p}\}$ be an IFS, where every contraction $f_{i}$ is a holomorphic map defined on a square $S$ (called the initial square of the IFS) such that $f_{i}(S) \Subset S$. The dynamical Cantor set $K$ associated to the IFS $\{f_{1},...,f_{p}\}$ is the limit set of the IFS $\{f_{1},...,f_{p}\}$, defined as follows :  
$$ K_{0} = S$$
$$  K_{n} = \bigcup_{1 \le i \le p} f_{i}(K_{n-1}) \text{ for every }n > 1$$  $$K =  \bigcap_{n \ge 0} K_{n}$$ \end{df57} We now introduce a formalism that will be used in the rest of the article : \begin{df57} Let $K$ be a dynamical Cantor set equal to the limit set of the IFS $\{f_{1},...,f_{p}\}$. A piece of depth $j$ is a connected component of $K_{j}$. Let $P$ be a  piece of depth $j$ and $Q$ a piece of depth $k$. We will say that :
\begin{enumerate}
  \item $P$ is the father of $Q$ if $k= j+1$ and if $Q  \Subset P$. We also say that $Q$ is a son of $P$  \item For every sequence $I = (i_{1},...,i_{k})$, we denote :$$K_{I} = f_{I}(K) = f_{i_{k}} \circ ... \circ f_{i_{1}} (K)$$ This is a dynamical Cantor set which is equal to the limit set of the IFS $\{f_{I} \circ f_{1} \circ f_{I}^{-1}, \ldots, f_{I} \circ f_{p}  \circ f_{I}^{-1}\}$ whose maps are all defined on $f_{I}(S)$ 
 \end{enumerate} \end{df57} In the following, for any IFS $\{f_{1},...,f_{p}\}$ of initial square $S$, we will suppose by simplicity that the center of $S$ is equal to 0.

  \begin{df57} Let  $\{f_{1},...,f_{p}\}$ be an IFS of initial square $S$. Let $P = f_{I}(S)$ be a piece of depth $k$ for some finite sequence $I$ of $k$ digits. \begin{enumerate}

\item 
The inscribed disk of $P$ is the disk of maximal diameter included in $P$ whose center is $f_{I}(0)$. We denote by $\delta(P)$ its diameter. \item  The middle inscribed disk is the disk of same center as the inscribed disk and of diameter multiplied by $\frac{1}{2}$. 
\item The escribed disk of $P$ is the disk of minimal diameter containing $P$ whose center is $f_{I}(0)$. We denote by $\Delta(P)$ its diameter.

\end{enumerate}

\end{df57}

\begin{rem57} These disks are well defined. For example, to define the diameter of the inscribed disk, it is enough to take the upper bound of the diameters of disks of center $f_{I}(0)$ included in $P$. Beware that for a square $S$ of diameter 1, we have $\delta(S) = \frac{1}{\sqrt{2}}$ and $\Delta(S) =1$. \end{rem57}

In the following we will define several notions associated to a fixed IFS $\{f_{1},...,f_{p}\}$. By simplification, they will be denoted relatively to the associated limit set, which will always be 
a dynamical Cantor set. For technical reasons, we will always suppose the existence of a disk $S'$ of center 0  of diameter larger than $\frac{1}{r}$ times the diameter of the square $S$ (for some $0<r<1$) such that each contraction $f_{i}$ can be extended in a contraction defined on $S'$ with $f_{i}(S') \Subset S$.

  \begin{df57} Let $f$ be a univalent map defined on $S'$ such that $f(S') \Subset S$. The distorsion of $f$ (on $S$), denoted by $D_{f}$, is defined by : 
 $$D_{f} =   \max_{(z,z') \in \overline{S}^{2}}   \left| \frac{f'(z)}{f'(z')} \right| $$
 We have : $1 \le D_{f}<+\infty$ according to the Koebe distorsion Theorem (since $S \Subset S'$). For any dynamical Cantor set $K$ associated to the IFS $\{f_{1}, \ldots,f_{p}\}$, we define the distorsion of $K$ to be the following number :

$$D_{K} = \sup_{I}  D_{f_{I}}  = \sup_{I} \Big(  \max_{(z,z') \in \overline{S}^{2}}  \left |\frac{ f'_{I}(z)}{f'_{I}(z')} \right| \Big) = \sup_{n \ge 0} \Big( \max_{|I| = n} \Big( \max_{(z,z') \in \overline{S}^{2}}  \left| \frac{ f'_{I}(z)}{f'_{I}(z')} \right| \Big)  \Big) $$ where the upper bound is taken on the set  of finite sequences of digits in $\{1, \ldots, p\}$.

\end{df57}

\begin{rem57}
The distorsion $D_{K}$ can be easily computed. Indeed, each contraction $f_{i}$ can be extended in a contraction defined on $S'$ with $f_{i}(S') \Subset S$ where $S'$ is a disk of center 0  of diameter larger than $\frac{1}{r}$ times the diameter of the square $S$ for some $0<r<1$. For every $I=(i_{1},\ldots,i_{k})$, we have $f_{I}(S') \subset f_{i_{k}}(S') \Subset S$. Then, by the Koebe distorsion Theorem we get  for every $I$ the uniform bounds : $$ \frac{1-r}{(1+r)^{3}} \le \left|\frac{ f'_{I}(z)}{f'_{I}(z')} \right| \le \frac{1+r}{(1-r)^{3}}$$ In particular, this implies $1 \le  D_{K} \le \min(    \frac{(1+r)^{3}}{1-r}      ,   \frac{1+r}{(1-r)^{3}})<+\infty$

\end{rem57} We will need the following technical lemma  : 
\begin{lem57} \label{torsi} 

For every finite sequence of digits $I$ in $\{1, \ldots, p\}$, we have that :
$$  \frac{1}{D_{K}} |f'_{I}(0)| \delta(S)  \le \delta(f_{I}(S)) \le D_{K} |f'_{I}(0)| \delta(S)$$
 $$  \frac{1}{D_{K}} |f'_{I}(0)| \Delta(S)  \le \Delta(f_{I}(S)) \le D_{K} |f'_{I}(0)| \Delta(S)$$
In particular, since $\frac{\Delta(S)}{\delta(S)} = \sqrt{2}$, we have : $\frac{1}{D_{K}^{2}}\sqrt{2} \le  \frac{\Delta(f_{I}(S))}{\delta(f_{I}(S))} \le D_{K}^{2} \sqrt{2}$

 \end{lem57} 
\begin{proof} The lemma is an easy consequence of the mean value inequality. The proof is left to the reader. 

  \end{proof}

The following lemma is an intermediate result : 
\begin{lem57} \label{circ} 
Let $z \neq z' \in S$. We denote by $\Gamma$ the disk of center $z$ such that $z' \in \partial \Gamma$. Then there exists a disk $\Gamma'$ of radius larger than $\frac{1}{5} |z-z'|$ included in $\Gamma \cap S$. 
\end{lem57}

\begin{proof} We can suppose by simplicity that the center of $S$ is 0 and that its diameter is equal to $\sqrt{2}$ (that is the length of a side of $S$ is 1 and $\delta(S) = 1$). Let us first suppose that $ |z-z'| \le \frac{1}{2}$. It is easy to check that this implies that $S \cap \Gamma$ contains at least one of the four quarters of the disk $\Gamma$. This quarter of disk contains a right triangle of sides $|z-z'|$, $|z-z'|$ and $\sqrt{2}|z-z'|$. The inscribed circle $\Gamma'$ of this triangle has its diameter equal to : $$2\frac{ 1 \times 1 } { 1 + 1+ \sqrt{2}} |z-z'| =\frac{2}{2+\sqrt{2}} |z-z'| $$ In the general case, $ |z-z'| \le \sqrt{2}$ and then $\Gamma$ contains a disk $\Gamma''$ of radius $\frac{1}{2\sqrt{2}}  |z-z'|  \le \frac{1}{2}$ of center $z$. This disk contains a disk $\Gamma'$ included in $S$ of diameter $\frac{1}{2\sqrt{2}} \frac{2}{2+\sqrt{2}} |z-z'|  = \frac{1}{2\sqrt{2}+2} |z-z'|> \frac{1}{5} |z-z'| $. Since $\Gamma'$ is included inside $\Gamma \cap S$, this concludes the proof.

\end{proof}

We will use the following lemma in the proof of Theorem A :

\begin{lem57} \label{hep}
Let $\Gamma_{1}$ be a disk of diameter 1. We denote by $\Gamma_{2}$ the disk of same center as $\Gamma_{1}$ and half diameter. Let $I$ be a finite sequence of digits in $\{1, \ldots, p\}$. We suppose that $f_{I}(S)$ intersects both $\Gamma_{2}$ and the complement of $\Gamma_{1}$. Then there exists a disk $\Gamma_{3}$ of diameter larger than $\frac{1}{20D_{K}^{2}}$ included in $\Gamma_{1}	 \cap f_{I}(S)$.

\end{lem57}

\begin{proof} 
Let us denote by $z_{1}$ a point of $\Gamma_{2} \cap f_{I}(S)$ and we take $z_{2} = f_{I}^{-1}(z_{1}) \in S$. We call $\mathcal{C} = \partial \Gamma_{1} \cap f_{I}(\overline{S})$ which is a compact set as an intersection of two compact sets. The set $\mathcal{C}' = f_{I}^{-1}(\mathcal{C})$ is then a compact set included in $\overline{S}$. We denote by $z_{3}$ a point of $\mathcal{C}'$ such that $|z_{3} -z_{2}| = \min_{z \in \mathcal{C}'} |z-z_{2}|$ (such a point $z_{3}$ exists by compacity). We denote by $\Omega_{1}$ the disk of center $z_{2}$ going through $z_{3}$. According to Lemma \ref{circ} (with $z_{2}$ as $z$, $z_{3}$ as $z'$, $\Omega_{1}$ as $\Gamma$ and $\Omega_{2}$ as $\Gamma'$), there exists a disk $\Omega_{2}$ of radius $\frac{1}{5} |z_{2}-z_{3}|$ included in $\Omega_{1} \cap S$. Since $\Omega_{2}$ is included in $\Omega_{1}$, for every $z \in \Omega_{2}$, we have $f_{I}(z) \in \Gamma_{1}$. According to Lemma \ref{torsi}, $f_{I}(\Omega_{2})$ contains a disk of diameter larger than $\frac{1}{D_{K}} |f'_{I}(0) | \cdot \frac{1}{5} |z_{2}-z_{3}|$. But according to Lemma \ref{torsi}, we also have : $ \frac{1}{4} \le  |f_{I}(z_{2})-f_{I}(z_{3})| \le D_{K} |f'_{I}(0) |  |z_{2}-z_{3}|$. Finally $f_{I}(\Omega_{2})$ contains a disk $\Gamma_{3}$ of diameter larger than $\frac{1}{20D_{K}^{2}}$ included in $\Gamma_{1}	 \cap f_{I}(S)$.

\end{proof}

  \begin{df57} \label{defff}  For every dynamical Cantor set $K$ of initial square $S$ associated to the IFS $\{f_{1},...,f_{p}\}$, we denote by $S_{i} = f_{i}(S)$ and we define the maximal reduction ratio $\lambda_{K}^{0}$ and the minimal reduction ratio $\Lambda_{K}^{0} $ of $K$ as follows :  $$\lambda_{K}^{0} =  \min_{i \in \{1,...,p\}} \frac {\delta(S_{i})} {\delta(S)}$$  
  $$\Lambda_{K}^{0} =  \max_{i \in \{1,...,p\}} \frac {\Delta(S_{i})} {\delta(S)}$$

 \end{df57}

 \begin{df57}
 For every dynamical Cantor set $K$ of initial domain $S$ associated to the IFS $\{f_{1},...,f_{p}\}$ of distorsion $D_{K}$, we define the distorted maximal reduction ratio $\lambda_{K}$ :  $$\lambda_{K} = \frac{1}{D_{K}^{2}} \cdot \min_{i \in \{1,...,p\}} \frac {\delta(S_{i})} {\delta(S)}$$ \end{df57} For every sequence $I$, we have $\lambda_{K_{I}}^{0}  \ge \lambda_{K}$ (note that $\lambda_{K_{I}}^{0}$ can be defined exactly as in Definition \ref{defff} even if the initial domain $f_{I}(S)$ is not a square). Informally speaking, in the construction of $K$, pieces are contracted by a factor larger than $\lambda_{K}$ from one step to another (from a father to a son). Indeed, by Lemma \ref{torsi} : $$\delta(f_{I}(S)) \le |f'_{I}(0)|D_{K} \delta(S)$$ and for every $i \in \{1, \ldots, p \}$ we have : $$\delta(f_{iI}(S)) \ge \frac{|f'_{I}(0)|}{D_{K}} \delta(S_{i})$$ 
$$\min_{1 \le i \le p}  \delta(f_{iI}(S)) \ge \frac{|f'_{I}(0)|}{D_{K}}  \min_{1 \le i \le p} \delta(S_{i})$$
Then : $$\lambda_{K_{I}}^{0} = \min_{i \in \{1,...,p\}} \frac {\delta (   f_{iI}(S)  )} {\delta(f_{I}(S))}      \ge  \frac{1}{D_{K}^{2}}  \min_{i \in \{1,...,p\}} \frac {\delta( S_{i}  )} {\delta ( S   )}     $$

 $$ \lambda_{K_{I}}^{0}  \ge \lambda_{K} $$  Similarly we can define : \begin{df57}
$$\Lambda_{K} = D_{K}^{2} \cdot \max_{i \in \{1,...,p\}} \frac {\Delta(S_{i})} {\delta(S)}$$ 
\end{df57} Informally speaking, in the construction of $K$, pieces are contracted by a factor smaller than $\Lambda_{K}$ from one step to another (from a father to a son). \medskip

We want to define a relevant notion of thickness. Similarly to the one-dimensional case, it will be defined as a quotient. The numerator is intended to measure the size of the pieces, this will be $\lambda_{K}$. The denominator will measure the size of the gap between points of the Cantor set, we define it now.

  \begin{df57}
 
 Let $K$ be a dynamical Cantor set $K$ of initial square $S$. We denote by $\rho(S) = 2 \max_{z \in S} \mathrm{dist}(z,K)$ the diameter of the largest disk included in $S$ which does not intersect $K$. We define the gap of $K$ :
                    $$\sigma^{0}_{K} = \frac{\rho(S)}{\delta(S)}$$ The Koebe distorsion Theorem implies that for every sequence $I$, we have  $\sigma^{0}_{    K_{I}}    \le D^{2}_{K}\sigma^{0}_{K}$ (again, $\sigma^{0}_{    K_{I}} $ can be defined with the same formula even if $f_{I}(S)$ is not a square). We set $\sigma_{K} = D_{K}^{2}\sigma^{0}_{K}$. 

\end{df57}

In particular, by definition, if we take any piece $P$ of $K$, every disk included in $P$ of diameter larger than $\sigma_{K}\delta(P)$ contains a point of $K \cap P$. \medskip

 Here is the definition of the thickness of a dynamical Cantor set : 

 \begin{df57} \label{deff}

Let $K$ be a dynamical Cantor set. We define the thickness $\mathrm{t}(K)$ of $K$ by :    $$\mathrm{t}(K) = \frac{\lambda_{K}}{D_{K}^{2}\sqrt{\sigma_{K}}}=\frac{\lambda_{K}^{0}}{D_{K}^{5}\sqrt{\sigma^{0}_{K}  }}  $$
\end{df57}

\begin{rem57}
We already saw that it is possible to give bounds for $D_{K}$. Since it is easy to estimate $\lambda_{K}^{0}$ and $\sigma^{0}_{K}$, it is also the case for $\mathrm{t}(K)$. 
\end{rem57}

We will work under the following condition : 

\begin{df57} Two dynamical Cantor sets $K$ and $L$ which are limit sets of respective IFS's defined on the same square are well balanced if : \begin{enumerate} \item  $\max(\Lambda_{K},\Lambda_{L})<\frac{1}{20}$ \item for every piece $P$ of $K$, every disk of diameter larger than $\frac{1}{20D_{L}^{2}}\delta(P)$ included in $P$ contains a son of $P$ \item for every piece $P$ of $L$, every disk of diameter larger than $\frac{1}{20D_{K}^{2}}\delta(P)$ included in $P$   contains a son of $P$

\end{enumerate}
\end{df57}

This condition simply means that : \begin{enumerate}  \item all the pieces are small compared to the initial square \item  the pieces occupy a substantial part of the initial square, there is not a big gap inside this initial square (up to a factor depending on the distorsion of the other Cantor set) \end{enumerate}

\begin{rem57}

A sufficient condition to satisfy condition 2. in the previous definition is $\sigma_{K}+2\Lambda_{K}< \frac{1}{20D_{L}^{2}}$ and similarly a sufficient condition to satisfy condition 3. is $\sigma_{L}+2\Lambda_{L}< \frac{1}{20D_{K}^{2}}$. Then a sufficient condition to satisfy the well balanced condition is : $$\max( \Lambda_{K},\Lambda_{L}, \sigma_{K}+2\Lambda_{K}, \sigma_{L}+2\Lambda_{L}) < \frac{1}{20\max(D_{K},D_{L})^{2}}$$
When the IFS is affine, which is often close to be the case, the distorsions $D_{K}$ and $D_{L}$ are equal to 1 : $\max( \Lambda_{K},\Lambda_{L}, \sigma_{K}+2\Lambda_{K}, \sigma_{L}+2\Lambda_{L}) < \frac{1}{20}$

\end{rem57}

\section{An example of robust intersection}

Before proving our results, let us give an example. In particular, it is intended to emphasize the fact that we are taking a square as an initial set in order to get easily arbitrarily large thicknesses with affine IFS's. We take the dynamical Cantor set $K$ which is equal to the limit set of the IFS $\{f_{1}, \ldots,f_{N^{2}}\}$ where every map $f_{k}$ ($1 \le k \le N^{2}$) is affine and  defined on a square $K_{0}$ of diameter 1. The set $K_{1} = \bigcup_{1 \le k \le N^{2}} f_{k}(K_{0})$ is the union of $N^{2}$ subsquares regularly located in the interior of $K_{0}$, of  relative size $0 < r < 1$. Then every subsquare has an inscribed disk of diameter equal to $\frac{1}{\sqrt{2}} \frac{r}{N}$. The limit set of this IFS is a dynamical Cantor set $K$. Since each $f_{k}$ is affine, we have $D_{K} = 1$. Since each subsquare has its inscribed disk of diameter equal to $\frac{1}{\sqrt{2}} \frac{r}{N}$ and the initial square has its inscribed disk of diameter $\frac{1}{\sqrt{2}}$, we have $\lambda_{K}^{0} =   \frac{r}{N}$. Every point of $K_{0}$ is distant from one of the subsquares at most $\frac{\sqrt{2}}{2}\frac{1-r}{N}$. In any of these subsquares, every point is distant at most  $\frac{\sqrt{2}}{2} \frac{1-r}{N} \cdot \frac{r}{N} $ from one of the subsubsquares and so on... Finally we have : $$\sigma_{K}^{0} = \frac{ \rho(K_{0})}{1/\sqrt{2}} \le  2\sqrt{2} \Big( \frac{\sqrt{2}}{2}  \frac{1-r}{N}+ \frac{\sqrt{2}}{2}  \frac{1-r}{N} \frac{r}{N}+  \frac{\sqrt{2}}{2}   \frac{1-r}{N} \Big( \frac{r}{N} \Big) ^{2} + \ldots  \Big) \le 4  \frac{1-r}{N}$$

$$\mathrm{t}(K) \ge \frac{r}{N \sqrt{ 4 (1-r)/N}}  = \frac{r}{\sqrt{ 4   (1-r)N}}$$ We take $N = 100 $ and $r =\frac{999}{1000}$. It is easy to check that $K$ and $K$ are well balanced and that we have $\mathrm{t}(K) > 1 $. In particular, we have $\mathrm{t}(K)^{2} > 1$. According to Corollary C, this implies that the two dynamical Cantor sets $K$ and $K$ have a robust intersection (i.e. if $f'_{k}$ and $f''_{k}$ are perturbations of $f_{k}$ then $K_{f'_{k}}$ and $K_{f''_{k}}$ intersect).

        \begin{center} \includegraphics[width = 10 cm]{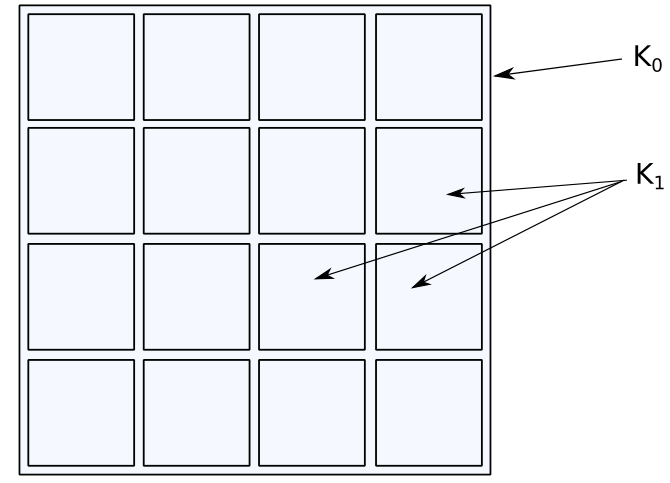} \end{center}  

$$\text{Figure 1 :  the sets $K_{0}$ and $K_{1}$ for $N=4$ and $r$ close to 1}  $$

\section{Proof of Theorem A}

We can suppose that the initial square is the square $S = S(0,1)$ centered at 0 of diameter 1. Indeed, the thickness of a dynamical Cantor set is defined from quotients of lengths and then invariant by rescaling, so it is possible to suppose that the initial square is $S(0,1)$. We also suppose that : $\sigma_{K} \le \sigma_{L}$. To show Theorem A, we construct by induction a sequence $(\alpha_{n})_{n \ge 0}$ of points of $K$ such that $\alpha_{n}  \in L_{n}$ for every $n \ge 0$, and more precisely $\alpha_{n}$ belongs to the middle inscribed disk of a piece of $L_{n}$   \medskip

We first show that this property is satisfied for $n = 0$. By hypothesis, $K$ and $L$ are well balanced so by definition there is a piece of $K$ of depth 1 included in every disk of diameter larger than $\frac{1}{20D^{2}_{L}}\delta(S)$ included in $S$. Since $\frac{1}{20D^{2}_{L}}\delta(S)\le \frac{1}{20} \frac{1}{\sqrt{2}}<\frac{1}{2\sqrt{2}}$, there exists a piece of depth 1 of $K$ included in the disk of center 0 of diameter $\frac{1}{2\sqrt{2}}$. But this disk is the middle inscribed disk of $S = S (0,1)$, which is the only piece of $L_{0}$. And since any piece of $K$ contains points of $K$, this implies that the middle inscribed disk of $S = S(0,1)$ contains a point $\alpha_{0}$ of $K$. Then the property is true for $n = 0$. \medskip

Let us now suppose that the property is true for some integer $n$ : there exists a point $\alpha_{n} \in K$ in the middle inscribed disk of a piece of $L_{n}$. We denote by $P$ this piece. We also denote $\rho(P) = 2 \max_{z \in P} \mathrm{dist}(z,L)$ the diameter of the largest disk included in $P$ which does not intersect $L$. The point $\alpha_{n}$ is the intersection of a sequence of pieces of $K$ whose first term is $S(0,1)$. We denote by $R_{1},\ldots,R_{q}$ the different sons of $P$. We distinguish two cases. \newline \newline  \underline{ Case 1, large gap :  $\rho(P) \ge\max_{1 \le j \le q} \delta(R_{j})$} \medskip \medskip

In this case, informally speaking, the sons $R_{1},\ldots,R_{q}$ of $P$ are sufficiently distant inside $P$ and not too close to one another. The point $\alpha_{n}$ is the intersection of a sequence of pieces of $K$ whose first term is $S(0,1)$. We choose in this sequence of pieces the piece $Q$ defined as follows : \begin{enumerate} \item if $(1+2\sqrt{2})D_{L}^{2} \rho(P) \ge \frac{1}{\sqrt{2}}$, we set $Q = S(0,1)$ \item if  $(1+2\sqrt{2})D_{L}^{2} \rho(P) <\frac{1}{\sqrt{2}}$, we define $Q$ as the last piece in the sequence such that $\delta(Q) \ge (1+2\sqrt{2})D_{L}^{2} \rho(P)$  \end{enumerate} Let us show that $Q$ contains a son $R$ of $P$. If $Q = S(0,1)$, this is obvious. If not, $Q$ has its inscribed disk $\Gamma$ of diameter larger than $(1+2\sqrt{2})D_{L}^{2} \rho(P)$. Since $\rho(P) \ge \max_{1 \le j \le q} \delta(R_{j})$, the diameter of the inscribed disk of every son of $P$ is smaller than $\rho(P)$. We then have two cases, depending on $\Gamma$ is (fully) included in $P$ or not. \medskip

 Let us first suppose that $\Gamma$ is included in $P$. Since $(1+2\sqrt{2})D_{L}^{2} \rho(P)> \rho(P)$, by definition of $\rho(P)$ the disk of same center as $\Gamma$ of diameter $\rho(P)$ contains a point of $L$ in its closure. This point of $L$ belongs to some son $R$ of $P$ and also to the escribed disk of $R$. This escribed disk is of diameter bounded from above by $\sqrt{2}D_{L}^{2} \rho(P)$ by Lemma \ref{torsi}. But $\Gamma $ is of diameter larger than $(1+2\sqrt{2})D_{L}^{2} \rho(P)$, then $\Gamma$ contains the escribed disk of $R$. Then $\Gamma$ contains $R$ and since $\Gamma \subset Q$ then $Q$ contains $R$ too. \medskip

Let us now suppose that $\Gamma$ is not included in $P$. We denote by $\Gamma'$ the inscribed disk of $P$ and by $\Gamma''$ its middle inscribed disk. Then $\alpha_{n} \in \Gamma''$ and $Q$ intersects both $\Gamma''$ and the complement of $\Gamma'$.
Then $Q \cap \Gamma'$ contains a disk $\Gamma'''$ of diameter larger than $\frac{1}{20 D_{K}^{2}}\delta(P)$ by Lemma \ref{hep} (with $\Gamma'$ as $\Gamma_{1}$, $\Gamma''$ as $\Gamma_{2}$ and $\Gamma'''$ as $\Gamma_{3}$). By hypothesis $K$ and $L$ are well balanced, so there is a son of $P$ which is included in every disk of diameter larger than $\frac{1}{20 D_{K}^{2}}\delta(P)$ included in $P$. In particular, there is a son $R$ of $P$ which is included in the disk $\Gamma'''$. Since $\Gamma''' \subset Q$, $Q$ contains $R$. \medskip

This shows that in any case $Q$ contains a son $R$ of $P$.

      \begin{center} \includegraphics[width = 11 cm]{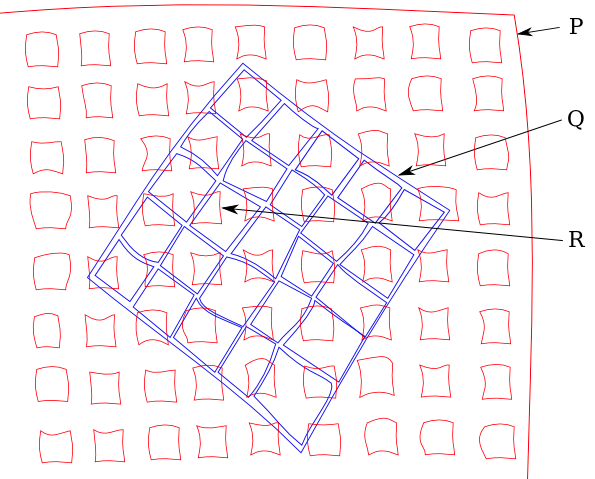} \end{center}  
$$\text{Figure 2 : construction of $R$ when $\Gamma \subset P$ } $$

We denote by $\rho(Q) = 2 \max_{z \in Q} \mathrm{dist}(z,K)$ the diameter of the largest disk included in $Q$ which does not intersect $K$. We have : $$   \frac{     \min_{U \text{ sons of }Q} \delta(U)        }{   \rho(Q)  }  \cdot \frac{     \min_{R_{j} \text{ sons of }P} \delta(R_{j})        }{ D_{L}^{2} \rho(P) } \ge \frac{    \lambda_{K}     }{    \sigma_{K} }  \cdot  \frac{    \lambda_{L}     }{ D_{L}^{2} \sigma_{L} }  \ge \frac{1}{ \sqrt{ \sigma_{K} \sigma_{L} }}   \cdot \frac{    \lambda_{K}     }{D_{K}^{2} \sqrt{\sigma_{K} }}  \cdot  \frac{    \lambda_{L}     }{D_{L}^{2} \sqrt{ \sigma_{L} }}  $$ where the last inequality comes from the fact that $D_{K} \ge 1$. Since $K$ and $L$ are well balanced, we have $\sigma_{K} < \frac{1}{20}$ and $\sigma_{L}< \frac{1}{20}$. We finally get : $$  \frac{     \min_{U \text{ sons of }Q} \delta(U)        }{   \rho(Q) }   \cdot \frac{     \min_{R_{j} \text{ sons of }P} \delta(R_{j})        }{D_{L}^{2}\rho(P) }  > 20 \cdot \text{t}(K)\text{t}(L)      \ge 20 \cdot 1 = 20$$ Since by definition of $Q$ we have $\delta(U_{k})<(1+2\sqrt{2})D_{L}^{2} \rho(P)$ for the son $U_{k}$ of $Q$ which contains $\alpha_{n}$, we get : $$ \min_{U \text{ sons of }Q} \delta(U)  <(1+2\sqrt{2})D_{L}^{2} \rho(P)$$  
$$\rho(Q) <   \frac{1+2\sqrt{2}}{20}   \min_{R_{j} \text{ sons of }P} \delta(R_{j})    <    \frac{1}{2}   \min_{R_{j} \text{ sons of }P} \delta(R_{j})  < \frac{1}{2}\delta(R) $$  Since $Q$ contains $R$, then there exists a point $\alpha_{n+1} \in K$ in the middle inscribed disk of the piece $R$ of $L_{n+1}$. The property is then satisfied for $n+1$.  \newline \newline \newline 
\underline{Case 2, small gap in $P$ : $\rho(P) < \max_{1 \le j \le q} \delta(R_{j})$} \medskip \medskip

In this case, informally speaking, the sons $R_{1},\ldots,R_{q}$ of $P$ are this time close to one another. For each son $R_{j}$ of $P$, we denote by $C_{j} = C_{j}^{x} \Subset R_{j}$ ($1 \le j \le q$) the disk of same center as the inscribed disk of $R_{j}$ and of diameter $\frac{1}{2}x$ times its diameter, where $0<x\le 1$ will be fixed in the following. For convenience, we will often write $C_{j}$ for $C_{j}^{x}$ when recalling the "$x$" is not needed but the disk $C_{j}$ depends on $x$.  \medskip 

We denote by $\tilde{\rho}(P)$ the diameter of the largest disk included in $P$ which does not intersect $C_{1} \cup \ldots \cup C_{q} = C_{1}^{x} \cup \ldots \cup C_{q}^{x}$. We have that  $\tilde{\rho}(P)$ is a continuous and decreasing function of $x$. We are now going to fix $0<x \le 1$. The disks $C_{1}^{1}, \ldots, C_{q}^{1}$ are the respective middle inscribed disks of $R_{1}, \ldots, R_{q}$. If we have $\max_{1 \le j \le q} \delta(C_{j}^{1}) \le \tilde{\rho}(P)$, then we fix $x = 1$. Let us suppose this is not the case. We have that $\max_{1 \le j \le q} \delta(C_{j}) = \max_{1 \le j \le q} \delta(C_{j}^{x})$ is a continuous and increasing function of $x$ which tends to 0 when $x$ tends to 0. On the other hand, $\tilde{\rho}(P)$ is a continuous function of $x$ strictly positive and decreasing. Then by the intermediate value Theorem there exists a value  $0<x<1$ such that $\max_{1 \le j \le q} \delta(C_{j}^{x}) = \tilde{\rho}(P)$ and we fix it. \medskip

 We want to give a bound on $\tilde{\rho}(P)$. Let us first suppose that $x = 1$. Let us take any disk $\Omega$ included in $P$ of diameter larger than $\max_{1 \le j \le q} \Delta(R_{j}) + \rho(P)$. Then the disk $\Omega'$ of diameter $\rho(P)$ of same center as $\Omega$ contains a point $\beta$ of $L$ by definition of $\rho(P)$. The point $\beta$ belongs to the escribed disk $\Omega''$ of a son $R_{l}$ of $P$. This escribed disk has its radius equal to $\frac{1}{2}\Delta(R_{l})$ with $\frac{1}{2}\Delta(R_{l}) \le \frac{1}{2} \max_{1 \le j \le q} \Delta(R_{j})$. Then the center of $\Omega''$ is at distance from $\beta$ at most $\frac{1}{2} \max_{1 \le j \le q} \Delta(R_{j})$. But the center of $\Omega''$ also belongs to $C_{l} = C_{l}^{1}$. Since the disk $\Omega$ has its diameter larger than $\max_{1 \le j \le q} \Delta(R_{j}) + \rho(P)$, $\Omega$ intersects $C_{1} \cup \ldots C_{q} = C^{1}_{1} \cup \ldots \cup C^{1}_{q}$. Then any disk included in $P$ of diameter larger than $\max_{1 \le j \le q} \Delta(R_{j}) + \rho(P)$ intersects $C_{1} \cup \ldots C_{q} = C^{1}_{1} \cup \ldots \cup C^{1}_{q}$. Then :  $$\tilde{\rho}(P) \le \max_{1 \le j \le q} \Delta(R_{j})+ \rho(P)$$According to Lemma \ref{torsi} and because $\rho(P) <  \max_{1 \le j \le q} \delta(R_{j})$, we have : $$\tilde{\rho}(P) \le  \sqrt{2}D_{L}^{2} \max_{1 \le j \le q} \delta(R_{j})+ \rho(P) < \sqrt{2}D_{L}^{2} \max_{1 \le j \le q} \delta(R_{j})+ \max_{1 \le j \le q} \delta(R_{j})$$ We have $D_{L} \ge 1$ and $\delta(R_{j}) = 2 \delta(C_{j}^{1})$ (remind that for $x =1$, $C^{1}_{j} = C_{j}$ is the middle inscribed disk of $R_{j}$). Then we get : $$\tilde{\rho}(P) < (1+\sqrt{2})D_{L}^{2} \max_{1 \le j \le q} \delta(R_{j}) =   2(1+\sqrt{2})D_{L}^{2} \max_{1 \le j \le q} \delta(C_{j})< 5 D_{L}^{2} \max_{1 \le j \le q} \delta(C_{j}) $$On the other hand, if $0<x<1$, we have $\max_{1 \le j \le q} \delta(C_{j}) = \tilde{\rho}(P)$. \medskip

No matter the value of $x$ we have : $$ \frac{\min_{1 \le j \le q} \delta(R_{j}) } { \max_{1 \le j \le q} \delta(R_{j}) } \ge \frac{\lambda_{L}}{\Lambda_{L}} >  20 \lambda_{L} $$  because $\Lambda_{L}<\frac{1}{20}$ (because $K$ and $L$ are well balanced). Then :  $$\min_{1 \le j \le q} \delta(C_{j}) = x \frac{1}{2} \min_{1 \le j \le q} \delta(R_{j}) > 20 \lambda_{L}   \cdot x \frac{1}{2} \max_{1 \le j \le q} \delta(R_{j}) = 20 \lambda_{L} \cdot
\max_{1 \le j \le q} \delta(C_{j}) $$ Then, no matter the value of $x$, we have :

 $$\frac{4}{ D_{L}^{2}} \lambda_{L} =   \frac{1}{5D_{L}^{2}} \cdot 20 \lambda_{L} = \frac{1}{\max(1,5D_{L}^{2})} \cdot 20 \lambda_{L}<  \frac{\min_{1 \le j \le q} \delta( C_{j} )}{\tilde{\rho}(P)}$$ From now on, the method will be the same as in Case 1, replacing $(P,R_{1}, \ldots, R_{q})$ by $(P,C_{1}, \ldots, C_{q})$. The point $\alpha_{n}$ is the intersection of a sequence of pieces of $K$ whose first term is $S(0,1)$. We choose in this sequence of pieces the piece $Q$ defined as follows :

\begin{enumerate} \item if $3\tilde{\rho}(P) \ge \frac{1}{\sqrt{2}}$, we set $Q = S(0,1)$ \item if  $3\tilde{\rho}(P) <\frac{1}{\sqrt{2}}$, we define $Q$ as the last piece in the sequence such that $\delta(Q) \ge 3\tilde{\rho}(P)$  \end{enumerate} Let us show that $Q$ contains one of the disks $C_{j}$. If $Q = S(0,1)$, this is obvious. If not, $Q$ has its inscribed disk $\Gamma$ of diameter larger than $3 \tilde{\rho}(P)$. Since $\tilde{\rho}(P) \ge \max_{1 \le j \le q} \delta(C_{j})$, the diameter of each of the disks $C_{j}$ is smaller than $\tilde{\rho}(P)$. We then have two cases, depending on $\Gamma$ is included in $P$ or not. \medskip

 Let us first suppose that $\Gamma$ is included in $P$. Since $3 \tilde{\rho}(P)> \tilde{\rho}(P)$, by definition of $\tilde{\rho}(P)$ the disk of same center as $\Gamma$ of diameter $\tilde{\rho}(P)$ contains a point of $C_{1} \cup \ldots \cup C_{q}$ in its closure. This point belongs to some disk $C_{k}$ which is of diameter bounded from above by $\tilde{\rho}(P)$. Since $\Gamma $ is of diameter larger than $ 3  \tilde{\rho}(P)$, then $\Gamma$ contains the disk $C_{k}$. Since $\Gamma \subset Q$ then $Q$ contains $C_{k}$ too. If $\Gamma$ is not included in $P$, the proof is similar as in Case 1 (pages 10 and 11) : $Q$ contains a son $R_{k}$ of $P$ and then also the disk $C_{k} \subset R_{k}$. \medskip

This shows that in any case $Q$ contains a son $R$ of $P$. We denote by $\rho(Q) = 2 \max_{z \in Q} \mathrm{dist}(z,K)$ the diameter of the largest disk included in $Q$ which does not intersect $K$. We have : 

$$    \frac{ \min_{U \text{ sons of }Q} \delta(U)  }{ \rho(Q) }  \cdot  \frac{\min_{1 \le j \le q} \delta(C_{j} )}{\tilde{\rho}(P)} \ge   \frac{ \lambda_{K} }{ \sigma_{K} }  \cdot \frac{4}{D_{L}^{2}} \lambda_{L} \ge  4 \frac{    \lambda_{K}     }{ \sqrt{ \sigma_{K}} }  \cdot   \frac{\lambda_{L}}{ D_{L}^{2} \sqrt{ \sigma_{L} }} $$ where the second inequality comes from the fact that $\sigma_{K} \le \sigma_{L}$ (this hypothesis was made at the beginning of the proof). Then we have : 
$$    \frac{     \min_{U \text{ sons of }Q} \delta(U)        }{ \rho(Q) } \cdot \frac{\min_{1 \le j \le q} \delta( C_{j} )}{\tilde{\rho}(P)} \ge  4 \frac{    \lambda_{K}     }{D_{K}^{2}\sqrt{ \sigma_{K}} } \cdot  \frac{\lambda_{L}}{D_{L}^{2}\sqrt{ \sigma_{L} }} \ge 4  \cdot  \mathrm{t}(K)\mathrm{t}(L)      \ge 4 $$ Since $\delta(U_{m})<3\tilde{\rho}(P)$ for the son $U_{m}$ of $Q$ which contains $\alpha_{n}$ (by definition of $Q$), we have :  $$\rho(Q)<  \frac{3}{4}\min_{1 \le j \le q} \delta( C_{j} ) < \min_{1 \le j \le q} \delta( C_{j} ) \le \delta(C_{k})$$  Since $C_{k}$ is included in $Q$, there exists a point $\alpha_{n+1} \in K$ in the disk $C_{k}$ and then in the middle inscribed disk of the corresponding $R_{k}$. The property is then true for $n+1$. By induction, it is true for every $n$. \medskip

In both cases we can conclude that $K \cap L_{n} \neq \emptyset$ for every integer $n$. This shows that $K \cap L$ is not empty. The proof of Theorem A is complete.

    \section{Proof of Theorem B}
Once defined the thickness $\mathrm{t}(K)$ of a dynamical Cantor set $K$, it is a natural question to investigate if it is a continuous function of the maps $f_{i}$ which define $K$. Indeed, as we already said it, Theorem A together with Theorem B give robust intersections of dynamical Cantor sets. We now prove Theorem B. To show this result, let us first remind that we have  : $ \mathrm{t}(K)  =  \frac{\lambda_{K}^{0}}{D_{K}^{5}\sqrt{\sigma^{0}_{K}  }}  $. To simplify, we will suppose that the initial domain of $K$ is the square $S =S(0,1)$ of center 0 of diameter 1. The following result is the key point in our proof :  \begin{lem57} \label{delta} The distorsion $D_{K}$ is a continuous function of the maps $f_{i}$. \end{lem57}
    
To show this, let us remind that : $$D_{K} = \sup_{I} \Big(  \max_{(z,z') \in \overline{S}^{2}} \frac{ |f'_{I}(z)|}{|f'_{I}(z')|} \Big) = \sup_{n \ge 0} \Big( \max_{|I| = n} \Big( \max_{(z,z') \in \overline{S}^{2}} \frac{ |f'_{I}(z)|}{|f'_{I}(z')|} \Big)  \Big) $$ In the following we will denote : $$S_{n}(f_{1},\ldots,f_{p}) = \max_{|I| = n} \Big( \max_{(z,z') \in \overline{S}^{2}} \frac{ |f'_{I}(z)|}{|f'_{I}(z')|} \Big) $$
We will need the following lemma : \begin{lem57} For every $\epsilon>0$, there exists an integer $n_{\epsilon} \ge 0 $ and a neighborhood $\mathcal{F}_{\epsilon}$ of $(f_{1}, \ldots, f_{p})$ of the form $\mathcal{F}_{\epsilon} = \{(g_{1}, \ldots, , g_{p}) \text{ such that }   : |g_{i}-f_{i}| <\eta_{\epsilon}\text{ on S } \}$ (for some constant $\eta_{\epsilon}$) such that : for every $n \ge n_{\epsilon}$, for every $(g_{1}, \ldots, g_{p}) \in \mathcal{F}_{\epsilon}$, we have :  $$  S_{n}(g_{1}, \ldots, g_{p})  \le S_{n_{\epsilon}}(g_{1}, \ldots, g_{p}) \cdot (1+\epsilon)$$  
\end{lem57} \begin{proof} We begin by fixing some constant $\eta_{\epsilon}>0$ and some neighborhood $\mathcal{F}_{\epsilon}= \{(g_{1}, \ldots, , g_{p}) \text{ such that }   : |g_{i}-f_{i}| <\eta_{\epsilon} \text{ on S } \}$ of $(f_{1}, \ldots, f_{p})$ and some constant $0<a<1$ such that for every $(g_{1}, \ldots, g_{p}) \in \mathcal{F}_{\epsilon}$, the IFS $\{ g_{1}, \ldots, g_{p} \}$ defines a dynamical Cantor set such that the diameter $\Delta(P)$ of the escribed disk of a piece $P$ is multiplied by no more than $a$ from one step to another. Such a constant $a$ exists for $(g_{1}, \ldots, g_{p})$ near  $(f_{1}, \ldots, f_{p})$ by contraction of the Poincar\'e metric, and it is easy to check that it can be taken locally constant. We also fix another constant :  reducing $\mathcal{F}_{\epsilon}$ if necessary, there exists $R>0$ such that for every $(g_{1}, \ldots, g_{p}) \in \mathcal{F}_{\epsilon}$, for every $z \in g_{1}(S) \cup \ldots \cup g_{p}(S)$, the disk of center $z$ of diameter $2R$ is included in $S$. Let us take any $(z,z') \in \overline{S}^{2}$. For every finite sequence $I =(i_{1},\ldots,i_{n})$, we have that : 

$$ \frac{ |g'_{I}(z)|}{|g'_{I}(z')|}  =  \frac{ |g'_{i_{1}}(z)|}{|g'_{i_{1}}(z')|} \cdot  \frac{ |g'_{i_{2}}(g_{i_{1}}(z))|}{|g'_{i_{2}}(g_{i_{1}}(z'))|} \cdot \cdots \cdot  \frac{ |g'_{i_{n}}(g_{i_{n-1}} \circ \ldots \circ g_{ i_{1}}(z))|}{|g'_{i_{n}}(g_{i_{n-1}} \circ \ldots \circ g_{ i_{1}}(z'))|}    $$ But for every $1 \le k \le n-1$, we have : $$  g_{i_{k}} \circ \ldots \circ g_{ i_{1}}(z) \in (g_{i_{k}} \circ \ldots \circ g_{i_{1}})(S) \text{ and } g_{i_{k}} \circ \ldots \circ g_{ i_{1}}(z') \in (g_{i_{k}} \circ \ldots \circ g_{i_{1}})(S)$$ and $\Delta((g_{i_{k}} \circ \ldots \circ g_{i_{1}})(S)) \le a^{k}$. Let us consider the restriction of $g_{i_{k+1}}$ to the ball of center $g_{i_{k}} \circ \ldots \circ g_{ i_{1}}(z)$ and of diameter $2R$ which is included in $S$ by  definition. The Koebe Theorem gives the following inequality : $$ \frac{ |g'_{i_{k+1}}(  g_{i_{k}} \circ \ldots \circ g_{ i_{1}}(z)   )|}{|g'_{i_{k+1}}(   g_{i_{k}} \circ \ldots \circ g_{ i_{1}}(z')  )|}  \le \frac{ 1 + a^{k}/R }  { (1- a^{k}/R)^{3}}  \le 1 + \frac{5}{R} a^{k} $$ if $k$ is large enough. The series $\sum_{k} a^{k} $ is absolutely converging, then the infinite product $\prod_{k} (1+\frac{5}{R} a^{k}) $ is also absolutely converging. In particular, when $n$ tends to $+ \infty$, $\prod_{k \ge n} (1+\frac{5}{R}a^{k}) $ tends to 1. It is then possible to take an integer $n_{\epsilon}$ such that $ \prod_{k \ge n_{\epsilon}} (1+\frac{5}{R}a^{k}) \le 1+\epsilon$. We fix this integer and then for every $n \ge n_{\epsilon}$ : $$ \frac{ |g'_{I}(z)|}{|g'_{I}(z')|}  =  \frac{ |g'_{i_{1}}(z)|}{|g'_{i_{1}}(z')|} \cdot  \frac{ |g'_{i_{2}}(g_{i_{1}}(z))|}{|g'_{i_{2}}(g_{i_{1}}(z'))|} \cdot \ldots \cdot  \frac{ |g'_{i_{n}}(     g_{i_{n-1}} \circ \ldots \circ g_{ i_{1}}(z)     )|}{|g'_{i_{n}}(   g_{i_{n-1}} \circ \ldots \circ g_{ i_{1}}(z')     )|}    $$

 $$ \frac{ |g'_{I}(z)|}{|g'_{I}(z')|}  \le  \frac{ |g'_{i_{1}}(z)|}{|g'_{i_{1}}(z')|} \cdot  \frac{ |g'_{i_{2}}(g_{i_{1}}(z))|}{|g'_{i_{2}}(g_{i_{1}}(z'))|} \cdot \ldots \cdot  \frac{ |g'_{i_{n_{\epsilon}}}(     g_{i_{n_{\epsilon}-1}} \circ \ldots \circ g_{ i_{1}}(z)             )|}{|g'_{i_{n_{\epsilon}}}(           g_{i_{n_{\epsilon}-1}} \circ \ldots \circ g_{ i_{1}}(z') )|}  \cdot (1+\epsilon)    $$  $$    \frac{ |g'_{I}(z)|}{|g'_{I}(z')|}  \le  S_{n_{\epsilon}}(g_{1}, \ldots, g_{p}) \cdot (1+\epsilon)$$ We take the upper bound on $(z,z') \in \overline{S}^{2}$ and on $|I| = n$ and we get : $$ S_{n}(g_{1}, \ldots, g_{p})  \le S_{n_{\epsilon}}(g_{1}, \ldots, g_{p}) \cdot (1+\epsilon) $$ The result follows.
\end{proof}

\begin{proof}[Proof of Lemma \ref{delta} : continuity of distorsion]  We now show that the distorsion $D_{K}$ is a continuous function of the maps $f_{i}$. Let us take any $\epsilon>0$. We begin by fixing the integer $n_{\epsilon}$ and the neighborhood $\mathcal{F}_{\epsilon}$ given by the previous lemma. In particular, for every $n \ge n_{\epsilon}$, for every $(g_{1}, \ldots, g_{p}) \in \mathcal{F}_{\epsilon}$, we have $  S_{n}(g_{1}, \ldots, g_{p})  \le S_{n_{\epsilon}}(g_{1}, \ldots, g_{p}) \cdot (1+\epsilon)$. Moreover, it is clear  that every term $S_{1}, \ldots, S_{n_{\epsilon}}$ is a continuous function of $(g_{1}, \ldots, g_{p}) \in \mathcal{F}_{\epsilon}$. It is then possible to choose a neighborhood $\mathcal{F}'_{\epsilon} \Subset \mathcal{F}_{\epsilon} $ such that for every $(g_{1}, \ldots, g_{p}) \in \mathcal{F}'_{\epsilon}$, we have : 
$$ \forall 1 \le k \le n_{\epsilon},  S_{k} (g_{1}, \ldots, g_{p})  \le D_{K} \cdot (1+\epsilon)$$
where $K$ is the dynamical Cantor set associated to the IFS $\{f_{1}, \ldots, f_{p}\} $. Let us denote by $L$ the dynamical Cantor set associated to the IFS $\{g_{1}, \ldots, g_{p}\} $. Then we have : $D_{L} \le D_{K} \cdot (1+\epsilon)^{2}$. \medskip

 We choose $l \ge 1$ such that $S_{l}(f_{1}, \ldots, f_{p})  \ge D_{K} \cdot (1-\epsilon)$. It is possible to take a new neighborhood $\mathcal{F}''_{\epsilon} \Subset \mathcal{F}'_{\epsilon} $ such that for every $(g_{1}, \ldots, g_{p}) \in \mathcal{F}''_{\epsilon}$, we have : $S_{l}(g_{1}, \ldots, g_{p})  \ge D_{K} \cdot (1-\epsilon)^{2}$ and then : $D_{L} \ge  D_{K} \cdot (1-\epsilon)^{2}$. \medskip

  As a conclusion, for every $(g_{1}, \ldots, g_{p}) \in \mathcal{F}''_{\epsilon}$, if $L$ is the limit set of $(g_{1}, \ldots, g_{p})$, we have :  $   D_{K} \cdot (1-\epsilon)^{2} \le  D_{L} \le  D_{K} \cdot (1+\epsilon)^{2}$.  This implies that the distorsion is a continuous function.

\end{proof}

\begin{proof}[Proof of Theorem B : continuity of the thickness] We saw that the distorsion $D_{K}$ is a continuous function of $f_{i}$. It is easy to check that $\lambda_{K}^{0}$ and $\sigma^{0}_{K}$ are also continuous functions of $f_{i}$, the proof is left to the reader. Finally, the thickness $\mathrm{t}(K)  =  \frac{\lambda_{K}^{0}}{D_{K}^{5}\sqrt{\sigma^{0}_{K}  }} $ is a continuous function of the contractions $f_{i}$.

\end{proof}

  \bibliographystyle{plain} \bibliography{bibi14}  \end{document}